\documentclass[12pt,twoside,reqno,a4paper]{amsart}

\usepackage[pdftex]{hyperref}
\usepackage{amsmath,amsthm,amsopn,amssymb,a4wide,times,varioref}
\usepackage{enumitem}
\usepackage[mathscr]{eucal}
\usepackage[utopia]{mathdesign}
\usepackage{tikz}
\usepackage[numbers,sort&compress]{natbib}
\newcommand{\ip}[2]{{
    \left<
      #1,#2
    \right>}}
\theoremstyle{proclaim}
\newtheorem{theorem}{Theorem}[section]
\newtheorem{lemma}[theorem]{Lemma}
\newtheorem{corollary}[theorem]{Corollary}
\newtheorem{proposition}[theorem]{Proposition}
\theoremstyle{statement}

\newtheorem{definition}[theorem]{Definition}
\newtheorem{example}[theorem]{Example}

\numberwithin{equation}{section}
\setlist{label={$($\roman{enumi}\kern1pt$)$}}
\DeclareFontFamily{U}{mathx}{\hyphenchar\font45}
\DeclareFontShape{U}{mathx}{m}{n}{
      <5> <6> <7> <8> <9> <10>
      <10.95> <12> <14.4> <17.28> <20.74> <24.88>
      mathx10
      }{}
\DeclareSymbolFont{mathx}{U}{mathx}{m}{n}
\DeclareFontSubstitution{U}{mathx}{m}{n}
\DeclareMathAccent{\widecheck}{0}{mathx}{"71}
\newcommand{\cla}{\mathcal{A}}
\newcommand{\clb}{\mathcal{B}}

\newcommand{\cll}{\mathcal{L}}

\title{Completely bounded kernels}
\author[Bhattacharyya]{Tirthankar Bhattacharyya${\vrule height8pt
    depth0pt width0pt}^*$}

\address{
  Department of Mathematics\\
  Indian Institute of Science\\
  Bangalore\\
  560 012\\
  India}

\email{tirtha@math.iisc.ernet.in}

\author[Dritschel]{Michael A.~Dritschel${\vrule height8pt
    depth0pt width0pt}^\dagger$}
\address{School of Mathematics {\&} Statistics\\
  Herschel Building\\ Newcastle University\\
  Newcastle upon Tyne\\
  NE1 7RU\\
  UK}

\email{m.a.dritschel@ncl.ac.uk}

\author[Todd]{Christopher S.~Todd}
\address{School of Mathematics {\&} Statistics\\
  Herschel Building\\ Newcastle University\\
  Newcastle upon Tyne\\
  NE1 7RU\\
  UK}

\email{c.s.todd@newcastle.ac.uk}

\thanks{${}^*$Research supported in part by grants from the UK-India
  Research and Education Initiative (UKIERI), DST (Ramanna Fellowship)
  and UGC SAP Phase IV}

\thanks{${}^\dagger$Research supported in part by grants
  from the UK-India Research and Education Initiative (UKIERI) and a
  Royal Society travel grant}

\dedicatory{Dedicated to the memory of Bela Sz.-Nagy, an inspiration
  to us all.}



\date{\today}

\begin{document}

\begin{abstract}
  We introduce completely bounded kernels taking values in $\cll (\cla
  , \clb)$ where $\cla$ and $\clb$ are $C^*$-algebras.  We show that
  if $\clb$ is injective such kernels have a Kolmogorov decomposition
  precisely when they can be scaled to be completely contractive, and
  that this is automatic when the index set is countable.
\end{abstract}

\maketitle

\setlength{\leftmargin}{0pt}

\section{Motivation}
\label{sec:motivation}

Decomposition properties of bounded maps play an important role in
functional analysis.  Some notable examples are
\begin{itemize}[label=$\bullet$]
\item The Hahn-Jordan decomposition of a complex measure in
    terms of a linear combination of four positive measures;
\item The decomposition of a bounded linear operator in terms of
  positive linear operators;
\item The Wittstock decomposition of a completely bounded linear map
  into an injective $C^*$-algebra in terms of a linear combination of
  four completely positive maps.
\end{itemize}

Fix a set $X$ and two \textbf{unital} $C^*$-algebras $\cla$ and
$\clb$.  Let $\cll ( \cla , \clb)$ be the space of bounded linear
maps from $\cla$ to $\clb$.  Let
\begin{equation*}
  {\mathbb K}_X ( \cla , \clb) = \{ k : X \times X \rightarrow \cll (
  \cla , \clb) \};
\end{equation*}
that is, $k$ is an $\cll ( \cla , \clb)$ valued kernel on $X$.  We
shall sometimes write ${\mathbb K}_X$ for $ {\mathbb K}_X ( \cla ,
\clb)$.

This set has an involution: if $k\in{\mathbb{K}_X(\mathcal{A,B})}$
then we define a kernel $k^{\ast}$ by
\begin{equation*}
  k^{\ast}(x,y)[a] = (k(y,x)[a^{\ast}])^{\ast}.
\end{equation*}
Notice that $(k^{\ast})^{\ast} = k$.  If $k=k^{\ast}$ then we call $k$
\textbf{hermitian} and denote the set of all hermitian $\cll ( \cla ,
\clb)$ valued kernels on $X$ by ${\mathbb K}^h_X = {\mathbb K}^h_X
(\cla , \clb)$.  Given $k \in {\mathbb K}_X$, the real and
imaginary parts of $k$ are two hermitian kernels in ${\mathbb K}_X$
defined by
\begin{equation*}
  \mathrm{Re}\,k = \tfrac{1}{2} (k + k^*) \quad\text{and}\quad
  \mathrm{Im}\, k = \tfrac{1}{2i} (k - k^*).
\end{equation*}
Thus $k =$ Re $k + i$ Im $k$.

The notation ${\mathbb K}^+$ will stand for \textbf{completely
  positive kernels}; that is, those $k \in \mathbb K$ which satisfy
\begin{equation}\label{def:cp_kernel}
  \sum_{i,j = 1}^n b_i^{\ast}k(x_i,x_j)[a_i^{\ast}a_j]b_j \geq {0}.
\end{equation}
for any $n \in \mathbb N$, any $x_1, x_2, \ldots x_n$ from $X$, $a_1,
a_2, \ldots a_n$ from $\cla$ and any $b_1, b_2, \ldots ,b_n$ from
$\clb$.  The definition originates in Barreto, Bhat, Liebscher and
Skeide \cite{BBLS03} and the name is justified by the equivalence of
(\textit{i\hspace{1pt}}) and (\textit{ii\hspace{1pt}}) in
Theorem~\ref{thm:positive_Kolmogorov} below.  Precursors are to be
found in papers by Murphy \cite{MU97}, where $C^*$-algebra valued
positive definite kernels are studied, and Heo \cite{Heo}, where the
set $X$ is assumed to be finite.  These papers make natural use of
Hilbert $C^*$-modules to give a so-called Kolmogorov decomposition of
a completely positive kernel.

\begin{theorem}[\cite{BBLS03}]
  \label{thm:positive_Kolmogorov}
  Let $k\in{\mathbb{K}}_X$.  Then the following are
  equivalent:
  \begin{enumerate}
  \item $k\in{\mathbb{K}}_X^{+}$.
  \item For any finite choice $x_1,x_2,\ldots x_n$ of elements from
    $X$ the (entrywise) map
    \begin{equation*}
      \left(k(x_i,x_j)\right) : M_n(\mathcal{A}) \to
      M_n(\mathcal{B})
    \end{equation*}
    is completely positive.
  \item The kernel $k$ has a \textbf{Kolmogorov decomposition}; that
    is there exists a pair $(E,\iota)$ with $E$ an
    $(\mathcal{A,B})$-correspondence and $\iota : X \to {E}$ a map
    such that for all choices of $x,y,a$,
    \begin{equation*}
      k(x,y)[a] = \ip{\iota(x)}{a\cdot\iota(y)}_{E}.
    \end{equation*}
  \end{enumerate}
\end{theorem}

As is standard in the literature for Hilbert $C^*$-modules, our inner
products are linear in the second argument.  A good recent exposition
on these concepts is the article by Skeide \cite{Skeide}.  For more
characterisations of completely positive kernels and a proof of
Theorem~\ref{thm:positive_Kolmogorov} see
\cite[Lemma~3.2.1,~Theorem~3.2.3]{BBLS03}.  We generalise the notion
of a Kolmogorov decomposition to non-positive kernels in
Section~\ref{sec:intr-cert-class}.

Let us write ${\mathbb D}^h_X = {\mathbb D}^h_X(\mathcal{A,B}) \subset
{\mathbb K}^h_X (\mathcal{A,B})$ for the set of hermitian
$\mathcal{L(A,B)}$ valued kernels on $X$ which can be expressed as the
difference of two completely positive kernels.  We call such hermitian
kernels \textbf{decomposable}.  More generally, a kernel in ${\mathbb
  K}_X (\mathcal{A,B})$ is said to be decomposable if its real and
imaginary parts are decomposable.  We shall see in
Theorem~\ref{thm:equiv_cond} that decomposability and the existence of
a form of Kolmogorov decomposition are equivalent, and so there is no
ambiguity in our use of the term ``decomposable''.  Note that the set
${\mathbb K}^h_X (\mathcal{A,B})$ carries a natural partial order: we
say $k_1 \leq {k_2}$ if $k_2 - k_1 \in {\mathbb{K}_X^{+}}$.

From part (\textit{iii\hspace{1pt}}) of
Theorem~\ref{thm:positive_Kolmogorov} it is clear that completely
positive kernels are hermitian.  Obviously, a difference of completely
positive kernels is also hermitian, and it is natural to wonder if
hermitian kernels can always be expressed as the difference of
completely positive kernels.  There are examples showing that this is
in general not possible.  For operator valued hermitian kernels, there
are characterisations of those with such a decomposition due to
Laurent Schwartz~\cite{SC64} which we generalize in
Theorem~\ref{thm:equiv_cond} (see also \cite{Alpay91,CG06}).  With
additional constraints all hermitian kernels can be decomposed.  For
example, \cite[Theorem~3.1]{Alpay91} implies that any $M_n(\mathbb C)$
valued hermitian kernel defined on a subset of $\mathbb C$ containing
a neighborhood of the origin on which it is jointly analytic is the
difference of positive kernels.  Similarly, by
\cite[Theorem~1.1.3]{ADRdS97}, an operator valued hermitian kernel
defined on an open subset of the complex plane with a finite number of
negative squares can be decomposed as the difference of positive
kernels.  This leads one to ask if there is some intrinsic property of
a kernel, hermitian or not, which corresponds to decomposability.  Our
goal in this paper is to determine precise conditions for the
decomposability of an $\mathcal{L(A,B)}$ valued kernel with $\clb$
injective.

All the examples at the very start of the paper can be viewed as
statements about kernels over one point sets.  More generally, kernels
appear quite naturally in many contexts both in analysis and its
applications, particularly in physics.  Decompositions of kernels are
key to the solution of a number of problems in which they appear, and
our own interest in the decomposition of non-positive kernels grew out
of the study of certain interpolation and realization problems.

Since completely positive maps are completely bounded, from
Theorem~\ref{thm:positive_Kolmogorov} a kernel which is decomposable
should be completely bounded when restricted to a finite subset of the
index set.  Also when the index set consists of one point the
Wittstock Theorem mentioned at the beginning equates decomposability
with complete boundedness.  We shall see that as a consequence of
Theorem~\ref{thm:equiv_cond}, that in general for kernels complete
boundedness (precise definition next section) is not quite sufficient,
and a sort of mild regularity condition, one which is automatically
satisfied by completely positive kernels and by completely bounded
kernels over countable index sets, must also be imposed.  The main
result of this paper is that complete boundedness plus regularity
characterises decomposable kernels, parallelling the Wittstock
theorems for completely bounded maps.  In the hermitian case this is
an appealing addition to the theorems of Schwartz, Constantinescu and
Gheondea.

We have organised the paper as follows.  In section~2, we introduce
the classes of kernels we shall be interested in, as well as our
generalisation of the notion of a Kolmogorov decomposition and the
definition of regularity mentioned above.  Section~3 generalises
results of Schwartz and Constantinescu and Gheondea while motivating
the regularity condition which decomposable kernels must satisfy.  It
also relates decomposability to a version of Paulsen's off-diagonal
technique, which is employed in Section~4 to give a series of
``local'' equivalences (i.e, those holding on finite index sets) and
``global'' equivalences (those holding on general index sets).
Section~5 then provides the topological tools needed to connect the
local and global results, yielding the main decomposability theorem in
Section~6.  We conclude in Section~7 with some examples and an
alternative, more local characterisation of the regularity condition.

\section{Introducing certain classes of kernels}
\label{sec:intr-cert-class}

In keeping with the theory of completely bounded maps, there is a
related notion of completely bounded kernel which we define as
follows.  Clearly, we are motivated by (\textit{ii\hspace{1pt}}) of
Theorem \ref{thm:positive_Kolmogorov}.

\begin{definition}
  \label{def:cb-kernel}
  Let $k : X\times{X} : \to \mathcal{L(A,B)}$.  If, given any finite
  subset $F = \{x_1,x_2,\ldots x_n\}$ of $X$, the (entrywise) map
  \begin{equation*}
    \left(k(x_i,x_j)\right) : M_n(\mathcal{A}) \to M_n(\mathcal{B})
  \end{equation*}
  is a completely bounded map, then we call $k$ a \textbf{completely
    bounded} kernel.
\end{definition}

This definition is not very restrictive as the following lemma
shows.

\begin{lemma}
  \label{lem:k-cb-iff-entries-cb}
  A kernel $k \in {\mathbb K}_X$ is completely bounded if and only if
  each $k(x,y)$ is a completely bounded map for $x,y \in X$.
\end{lemma}

\begin{proof}
  Clearly if $k$ is a completely bounded kernel, then each $k(x,y)$
  is a completely bounded map.  Conversely, if each $k(x,y)$ is a
  completely bounded map, then for any finite subset $F =
  \{x_1,x_2,\ldots x_n\}$ of $X$, the map
  \begin{equation*}
    \left(k(x_i,x_j)\right) : M_n(\mathcal{A}) \to M_n(\mathcal{B})
  \end{equation*}
  is a sum of completely bounded maps and hence is completely
  bounded.
\end{proof}

\begin{definition}
  \label{def:Kolmogorov-decomp}
  A kernel $k\in\mathcal{\mathbb{K}(A,B)}$ has a \textbf{Kolmogorov
    decomposition} if there exists a triple $(E,\iota,J)$ where $E$ is
  an $(\mathcal{A,B})$-correspondence, $\iota : X \to {E}$ a map, and
  $J\in\mathcal{L}^a(E)$ a contractive $\mathcal{A}$-module map, such
  that
  \begin{equation*}
    k(x,y)[a] = \ip{J \iota(x)}{ a\cdot \iota(y)}_{E} \qquad 
    \text{ for all } x \in X, \; y\in Y \text{ and } a \in \cla.
  \end{equation*}
  If $J$ is self-adjoint (respectively positive), we shall say $k$ has
  a \textbf{hermitian} (respectively \textbf{positive}) Kolmogorov
  decomposition.  A Kolmogorov decomposition $(E, \iota, J)$ of a
  kernel $k$ is said to be \textbf{minimal} if $E =
  \overline{\mathrm{span}}\, \cla \cdot \iota(X) \clb$.
\end{definition}

Theorem~\ref{thm:positive_Kolmogorov} states that a completely
positive kernel has a Kolmogorov decomposition with $J$ equal to the
identity map.  In general, if $k$ has a Kolmogorov decomposition, then
\begin{equation*}
  k^{\ast}(x,y)[a] = (k(y,x)[a^{\ast}])^{\ast} =
  \ip{J \iota(y)}{a^{\ast} \cdot \iota(x)}^{\ast} = 
  \ip{a^{\ast} \cdot \iota(x)}{J \iota(y)} =
  \ip{J^\ast \iota(x)}{a \cdot \iota(y)}.
\end{equation*}
Thus a kernel with a hermitian Kolmogorov decomposition is hermitian.
On the other hand, if $k$ is hermitian and $(E, \iota, J)$ is a
minimal Kolmogorov decomposition of $k$, then it is not difficult to
see that the Kolmogorov decomposition is hermitian (see
Theorem~\ref{thm:equiv_cond} below).

Given two minimal decompositions $(E, \iota, J)$ and $(E^\prime,
\iota^\prime, J^\prime)$, let $U : E \to E^\prime$ be defined by $U
\iota(x) = \iota^\prime(x)$ and extend bilinearly.  Then $U$
is a surjective isometry.  Also,
\begin{equation*}
  \ip{J \iota(x)}{a \cdot \iota(y)} = \ip{J^\prime \iota^\prime(x)}{a
    \cdot \iota^\prime(y)} = \ip{J^\prime U ( a \cdot \iota(x))}{U(a
    \cdot \iota(y))} = \ip{U^* J^\prime U ( a \cdot \iota(x))}{a \cdot
    \iota(y)},
\end{equation*}
hence $U^* J^\prime U = J$.  Completely positive kernels have minimal
Kolmogorov decompositions~\cite[Proposition~3.1.3]{BBLS03}.

\begin{definition}
  \label{def:unif-cb-and-cc}
  The kernel $k : X\times{X} : \to \mathcal{L(A,B)}$ is called a
  \textbf{uniformly completely bounded} kernel if there is a constant
  $M \geq 0$ such that for any finite subset $F = \{x_1,x_2,\ldots
  x_n\}$ of $X$, the map
  \begin{equation*}
    \left(k(x_i,x_j)\right) : M_n(\mathcal{A}) \to M_n(\mathcal{B})
  \end{equation*}
  is a completely bounded map with cb-norm at most $M$.  We say that
  $k$ is \textbf{completely contractive} if we can choose $M=1$.
\end{definition}

A completely positive kernel need not be uniformly completely bounded.
However, it can be rescaled to be so.  Indeed, a completely positive
kernel $k$ can be conjugated to give a completely contractive kernel.
If $k \in {\mathbb K}_X^+ (\cla, \clb)$ then there is a scalar valued
map $c:X\to (0,1]$ on $X$ such that the new kernel $\tilde{k}$ in
${\mathbb K}_X (\cla, \clb)$ defined by $\tilde{k}(x,y)[a]= c(x)
k(x,y)[a] \overline{c(y)}$ is a completely positive and completely
contractive kernel.  To see this, let the minimal positive Kolmogorov
decomposition of $k$ be $(E, \iota)$.  Define $c(x) =
\|\ip{i(x)}{i(x)} \|^{-1/2}$ or $1$ if $\|\ip{\iota(x)}{\iota(x)}\| <
1$.  Then for any natural number $n$ and any $x_1, x_2, \ldots ,x_n$
from $X$, the map
\begin{equation*}
  \left(\tilde{k}(x_i,x_j)\right) : M_n(\mathcal{A}) \to
  M_n(\mathcal{B})
\end{equation*}
is a completely positive map since
\begin{equation*}
  \left( \tilde{k} (x_i,x_j) \right) = 
  \mathrm{diag}\, \{ c(x_1), \ldots ,c(x_n) \}
  \left( k (x_i,x_j) \right) 
  \mathrm{diag}\, \{c(x_1), \ldots ,c(x_n) \}.
\end{equation*}
Also
\begin{equation*}
  \left\|\left( \tilde{k} (x_i,x_j) \right) I_{M_n{(\cla)}}\right\| = 
  \left\| \mathrm{diag}\, \left\{ c(x_1)^2\ip{i(x_1)}{i(x_1)}, 
      \ldots ,c(x_n)^2\ip{i(x_n)}{i(x_n)} \right\} \right\| \leq 1,
\end{equation*}
so $\tilde{k} $ is completely contractive.  Moreover, the minimal
positive Kolmogorov decomposition of $\tilde{k}$ is $(E,
\tilde{\iota})$ where $\tilde{\iota}(x)$ is a multiple of $\iota(x)$
for every $x \in X$.

This leads to the following definition.

\begin{definition}
  \label{def:regular-cb-kernel}
  A completely bounded kernel $k \in {\mathbb K}_X (\cla, \clb)$ is
  \textbf{regular} if there is a scalar valued map $c:X\to (0,1]$ such
  that $\tilde{k}(x,y)[a] := c(x) k(x,y)[a] \overline{c(y)}$ is
  completely contractive.
\end{definition}

As noted, completely positive kernels are automatically regular, and
as we shall see, regular completely bounded kernels are precisely the
ones with Kolmogorov decompositions.  In section~7 we give an example
of a completely bounded kernel which s not regular, as well as give
several instances where regularity is automatic.

\section{Kolmogorov Decomposition and Matrix Completion}
\label{sec:kolmogorov_decompositions}

The chief result of this section is Lemma
\ref{lem:kolm_decomp_characterisation} which will be used in the
proof of the main theorem of this paper in Section 6.

\begin{theorem}
  \label{thm:equiv_cond}
  Let $k\in {\mathbb{K}}_X$.  Then the following are equivalent:
  \begin{enumerate}
  \item $k\in{\mathbb D}_X$.
  \item There is a kernel $L \in {\mathbb K}_X^+$ such that $-L \leq k
    \leq {L}$.
    \item There is a kernel $L \in {\mathbb K}_X^+$ such that $k \leq
      {L}$.
  \item $k$ has a hermitian Kolmogorov decomposition.
  \item $k$ has a hermitian Kolmogorov decomposition $(E, \iota , J)$
    such that $J^2 = I_E$.
  \end{enumerate}
\end{theorem}

\begin{proof}
  (\emph{i}\kern1pt) $\Rightarrow$ (\emph{ii}\kern1pt): Let $k = k_1 -
  k_2$.  Then
  \begin{equation*}
    k + k_1 + k_2 = k + k_1 + (k_1 - k) = 2k_1 \geq {0}.
  \end{equation*}
  Thus $k \geq {-(k_1 + k_2)}$.  Similarly
  \begin{equation*}
    (k_1 + k_2) - k = (k + k_2) + k_2 - k = 2k_2 \geq {0}
  \end{equation*}
  so $k_1 + k_2 \geq {k}$.  Since $k_1 + k_2$ is completely positive,
  set $L = k_1 + k_2$ to obtain the result.

  (\emph{ii}\kern1pt) $\Rightarrow$ (\emph{iii}\kern1pt): Immediate.

  (\emph{iii}\kern1pt) $\Rightarrow$ (\emph{i}\kern1pt): Let $k \leq
  {L}$.  Then $k = L - (L - k)$ is the difference of completely
  positive kernels.

  (\emph{i}\kern1pt) $\Rightarrow$ (\emph{v}\kern1pt): This is a
  construction analogous to constructing decompositions of completely
  positive kernels.  Let $k = k_1 - k_2$ and let $(E_1, \iota_1, J_1)$
  and $(E_2 , \iota_2, J_2)$ be the Kolmogorov decompositions of the
  completely positive kernels $k_1$ and $k_2$ respectively.  Define
  $E: E_1\oplus{E_2}$ and $\iota(x) = \iota_1 (x) \oplus \iota_2 (x)$
  and $J : E \to {E}$ by $J(e_1\oplus{e_2}) = e_1\oplus (-e_2)$.  Then
  \begin{equation*}
    \begin{split}
      k(x,y)[a] = k_1(x,y)[a] - k_2(x,y)[a]&
      = \ip{ \iota_1(x)}{a\cdot \iota_1(y)} 
      - \ip{ \iota_2(x)}{a\cdot \iota_2(y)}\\
      & = \ip{J ( \iota_1(x) \oplus  \iota_2(x))}{a\cdot
        ( \iota_1(y) \oplus \iota_2(y))}\\
      & = \ip{J \iota(x) }{a\cdot \iota(y)}.
    \end{split}
  \end{equation*}
  Notice that $J$ is self-adjoint, a left $\mathcal{A}$-module map,
  and $J^2 = I$.  This completes the construction of the
  decomposition of $k$.

  (\emph{v}\kern1pt) $\Rightarrow$ (\emph{iv}\kern1pt): Immediate.

  (\emph{iv}\kern1pt) $\Rightarrow$ (\emph{i}\kern1pt): Since $J$ is a
  self-adjoint element of the $C^{\ast}$-algebra $\mathcal{L}^a(E)$
  it can be expressed as the difference of two positive elements, say
  $J = J_1 - J_2$.  Since $J$ is an $\mathcal{A}$-module map, each
  of $J_1,J_2$ must also be so.  Then we have
  \begin{equation*}
    \begin{split}
      k(x,y)[a]& = \ip{J \iota(x)}{a\cdot  \iota(y)}\\
      & = \ip{J_1 \iota(x) - J_2 \iota(x)}{a\cdot \iota(y) }\\
      & = \ip{J_1 \iota(x)}{a \cdot \iota(y)} 
      - \ip{J_2 \iota(x)}{a \cdot \iota(y)}\\
      & = \ip{ J_1^{1/2}  \iota(x) }{J_1^{1/2} ( a\cdot \iota(y))} -
      \ip{ J_2^{1/2} \iota(x) }{J_2^{1/2} ( a\cdot \iota(y))}\\
      & = \ip{ J_1^{1/2}  \iota(x) }{a\cdot (J_1^{1/2} \iota(y))} -
      \ip{ J_2^{1/2} \iota(x) }{a\cdot (J_2^{1/2}  \iota(y))}\\
    \end{split}
  \end{equation*}
  which gives decompositions of two completely positive kernels whose
  difference is $k$.
\end{proof}

 We now know that:
\begin{itemize}[label=$\bullet$]
\item A kernel $k$ is completely positive if and only if it has a
  positive Kolmogorov decomposition since $\ip{J \iota(x)}{a \cdot
    \iota(y)} = \ip{ J^{1/2}\iota(x)}{a \cdot J^{1/2}\iota(y)} = \ip{
    \widetilde{\iota}(x)}{a \cdot \widetilde{\iota}(y)}$).
\item A kernel $k$ is the difference of two completely positive
  kernels if and only if it has a hermitian Kolmogorov decomposition.
  This was the content of Theorem~\ref{thm:equiv_cond}.
\item If $k$ has a hermitian Kolmogorov decomposition then we can
  assume that the operator $J$ is unitary, since by
  Theorem~\ref{thm:equiv_cond} (\emph{v}\kern1pt) we can take $J^2 =
  I$.
\item Splitting a general kernel into real and imaginary parts, it is
  clear that such a kernel is a linear combination of at most four
  completely positive kernels if and only if it has a Kolmogorov
  decomposition.
\item Since completely positive kernels are regular, from the
  equivalence of (\emph{i}\kern1pt) and (\emph{ii}\kern1pt) in
  Theorem~\ref{thm:equiv_cond} we gather that decomposable hermitian
  kernels (and hence general decomposable kernels) are regular (this
  also follows easily from Corollary~\ref{cor:Haagerup_says_so}
  below).
\end{itemize}
We make use of these facts in proving the next lemma, a generalisation
of \cite[Theorem 4.4]{CG06}.  Compare also with Theorem~8.3
of \cite{PA02}.

\begin{lemma}
  \label{lem:kolm_decomp_characterisation}
  Let $k\in{\mathbb{K}_X(\mathcal{A,B})}$.  Then $k$ has a Kolmogorov
  decomposition if and only if there exist completely positive kernels
  $L_1$ and $L_2$ in ${\mathbb{K}_X(\mathcal{A,B})}$ such that
  \begin{equation*}
    (x,y) \mapsto
    \begin{pmatrix}
      L_1(x,y) & k(x,y)\\
      k^{\ast}(x,y) & L_2(x,y)
    \end{pmatrix}
    : \mathcal{A} \to M_2(\mathcal{B})
  \end{equation*}
  is a completely positive kernel.
\end{lemma}

\begin{proof}
  If $k(x,y)[a] = \ip{J \iota(x)}{a\cdot \iota(y)}$ then we take
  \begin{equation*}
    L_1(x,y)[a] = L_2(x,y)[a] = \ip{\iota(x)}{a\cdot \iota(y)}
  \end{equation*}
  which gives us
  \begin{equation*}
    \begin{pmatrix}
      L_1(x,y)[a] & k(x,y)[a]\\
      k^{\ast}(x,y)[a] & L_2(x,y)[a]
    \end{pmatrix}
    =
    \begin{pmatrix}
      \ip{\iota(x)}{a\cdot \iota(y)} & \ip{J \iota(x)}{a\cdot
        \iota(y)}
      \\[5pt]
      \ip{J^* \iota(x)}{a\cdot \iota(y)}
      & \ip{ \iota(x)}{a\cdot \iota(y)}
    \end{pmatrix}
  \end{equation*}
  View $M_2(E)$ as an $M_2(\mathcal{B})$-module.  The left
  $\mathcal{A}$-action is defined by embedding $\mathcal{A}$ in
  $\mathcal{L}^a(E)$, which can in turn be identified with the
  diagonal of $M_2(\mathcal{L}^a(E))$, which is completely
  isometrically isomorphic to $\mathcal{L}^a(M_2(E))$.  Notice that
  \begin{equation*}
    \begin{pmatrix}
      I & J\\
      J^* & I
    \end{pmatrix}
  \end{equation*}
  is then a positive element of $\mathcal{L}^a(M_2(E))$ and commutes
  with the left $\mathcal{A}$-action.  Combining this with the above
  gives us
  \begin{equation*}
    \begin{pmatrix}
      L_1(x,y)[a] & k(x,y)[a]\\
      k^{\ast}(x,y)[a] & L_2(x,y)[a]
    \end{pmatrix}
    = \ip{\begin{pmatrix}
        I & J \\
        J^* & I
      \end{pmatrix}
      \begin{pmatrix}
        \iota(x) & 0\\
        0 & \iota(x)
      \end{pmatrix}
    }{
a \cdot
      \begin{pmatrix}
        \iota(y) & 0\\
        0 & \iota(y)
      \end{pmatrix}
    }
  \end{equation*}
  which is a positive Kolmogorov decomposition.  Thus
  \begin{equation*}
    \begin{pmatrix}
      L_1 & k \\
      k^{\ast} & L_2
    \end{pmatrix}
  \end{equation*}
  is a completely positive kernel.

  Conversely suppose that the matrix of kernels is completely
  positive.  Then conjugation by an element of $M_2(\mathcal{B})$
  preserves complete positivity.  In particular
  \begin{equation*}
    {\begin{pmatrix}
        1 & 1\\
        1 & -1
      \end{pmatrix}}^{\ast}
    \begin{pmatrix}
      L_1 & k\\
      k^{\ast} & L_2
    \end{pmatrix}
    \begin{pmatrix}
      1 & 1\\
      1 & -1
    \end{pmatrix}
    =
    \begin{pmatrix}
      L_1 + L_2 + k + k^{\ast} & L_1 - L_2 - (k - k^{\ast})\\
      L_1 - L_2 + (k - k^{\ast}) & L_1 + L_2 - (k + k^{\ast})
    \end{pmatrix}
  \end{equation*}
  is a completely positive kernel.  The entries on the diagonal must
  be completely positive, from which we deduce that
  \begin{equation*}
    -\tfrac{1}{2}(L_1 + L_2) \leq \tfrac{1}{2}(k + k^{\ast}) \leq
    \tfrac{1}{2}(L_1 + L_2).
  \end{equation*}
  Similarly, conjugation by the $M_2(\mathcal{B})$ element
  \begin{equation*}
    \begin{pmatrix}
      1 & i\\
      -i & -1
    \end{pmatrix}
  \end{equation*}
  tells us that
  \begin{equation*}
    -\tfrac{1}{2}(L_1 + L_2) \leq \tfrac{-i}{2}(k - k^{\ast}) \leq
    \tfrac{1}{2}(L_1 + L_2).
  \end{equation*}
  By Theorem~\ref{thm:equiv_cond}, the kernels $K_1 = \tfrac{1}{2}(k +
  k^{\ast})$ and $K_2 = \tfrac{-i}{2}(k - k^{\ast})$ have hermitian
  Kolmogorov decompositions.  Let $(E_j, \iota_j, J_j)$ be the Kolmogorov
  decomposition of $K_j$ for $j-1,2$.  Now
  \begin{equation*}
    \begin{split}
      k(x,y)[a]& = K_1(x,y)[a] + i K_2(x,y)[a] \\
      & = \ip{J_1 {\iota_1(x)}}{a\cdot {\iota_1(y)}} +
      i\ip{J_2 {\iota_2(x)}}{a\cdot {\iota_2(y)}}\\
      & = \ip{(J_1 \oplus {(-iJ_2)}) (\iota_1(x)\oplus
      {\iota_2(x)})}{a \cdot
        (\iota_1(y)\oplus{\iota_2(y)})}_{E_1\oplus
        {E_2}}.
    \end{split}
  \end{equation*}
  Thus $k$ has a Kolmogorov decomposition.
\end{proof}

If $\varphi_1,\varphi_2,\varphi:\mathcal{A}\to\mathcal{B}$ and
we define
\begin{equation*}
  \Phi : M_2(\mathcal{A}) \to M_2(\mathcal{B}):
  \begin{pmatrix}
    a & b \\ c & d
  \end{pmatrix}
  \mapsto
  \begin{pmatrix}
    \varphi_1(a) & \varphi(b) \\ \varphi^{\ast}(c) & \varphi_2(d)
  \end{pmatrix}
  \qquad
  \Psi : \mathcal{A} \to M_2(\mathcal{B}) : a \mapsto
  \begin{pmatrix}
    \varphi_1(a)&\varphi(a)\\\varphi^{\ast}(a)&\varphi_2(a)
  \end{pmatrix}
\end{equation*}
then it is a result due to Haagerup that $\Phi$ is completely positive
if and only if $\Psi$ is completely positive~\cite{HA85} (though our
presentation follows~\cite{PA02}).  As a consequence, if we then allow
a $2\times{2}$ matrix of kernels to act (at each point) as a Schur
product
\begin{equation*}
  \begin{pmatrix}
    k_1 & k_2 \\ k_3 & k_4
  \end{pmatrix}
  (x,y)
  \left[
    \begin{pmatrix}
      a_{1,1} & a_{1,2} \\ a_{2,1} & a_{2,2}
    \end{pmatrix}
  \right]
  =
  \begin{pmatrix}
    k_1(x,y)[a_{1,1}] & k_2(x,y)[a_{1,2}] \\ k_3(x,y)[a_{2,1}] &
    k_4(x,y)[a_{2,2}]
  \end{pmatrix}
\end{equation*}
we obtain the following:

\begin{corollary}
  \label{cor:Haagerup_says_so}
  The kernel $k\in{\mathbb{K}(\mathcal{A,B})}$ has a Kolmogorov
  decomposition if and only if there exist completely positive kernels
  $L_1, L_2\in{\mathbb{K}^+(\mathcal{A,B})}$ such that
  \begin{equation*}
    (x,y) \mapsto
    \begin{pmatrix}
      L_1(x,y) & k(x,y)\\
      k^{\ast}(x,y) & L_2(x,y)
    \end{pmatrix}
    :M_2(\mathcal{A}) \to M_2(\mathcal{B})
  \end{equation*}
  is a completely positive kernel.
\end{corollary}

\section{An application of the off-diagonal method}
\label{sec:an-appl-diag}

In light of Corollary~\ref{cor:Haagerup_says_so}, the study of the
decomposability of a completely bounded map $\varphi$ is related to
the problem of completely positive completion of a $2\times 2$ matrix
with $\varphi$ and $\varphi^{\ast}$ in the off-diagonal positions.
Our goal is to use this relationship to show that any uniformly
completely bounded kernel (and hence any regular completely bounded
kernel) into an appropriate space has a Kolmogorov decomposition.
Now, let $k\in{\mathbb{K}_X(\mathcal{A,B})}$ and consider the
following six statements:

\begin{itemize}
\item[(\emph{i}\kern1pt)] There exist
  $L_1,L_2\in{\mathbb{K}^{+}_X(\mathcal{A,B})}$ such that
  \begin{equation*}
    \begin{pmatrix}
      L_1 & k\\
      k^{\ast} & L_2
    \end{pmatrix}
    \in \mathbb{K}^{+}_X(M_2(\mathcal{A}),M_2(\mathcal{B})).
  \end{equation*}
\item[(\emph{ii}\kern1pt)] There exist
  $L_1,L_2\in{\mathbb{K}^{+}_X(\mathcal{A,B})}$ such that given any
  finite subset $x_1,x_2,\ldots x_n$ of $X$ the map
  \begin{equation*}
    {\begin{pmatrix}
        L_1(x_i,x_j) & k(x_i,x_j)\\
        k^{\ast}(x_i,x_j) & L_2(x_i,x_j)
      \end{pmatrix}}_{i,j = 1}^n :
    M_n(M_2(\mathcal{A})) \to M_n(M_2(\mathcal{B}))
  \end{equation*}
  is completely positive.
\item[(\emph{iii}\kern1pt)] There exist
  $L_1,L_2\in{\mathbb{K}^{+}_X(\mathcal{A,B})}$ such that given any
  finite subset $x_1,x_2,\ldots x_n$ of $X$ the map
  \begin{equation*}
    \begin{pmatrix}
      \left(L_1(x_i,x_j)\right)_{i,j = 1}^n&
      \left(k(x_i,x_j)\right)_{i,j = 1}^n\\[6pt]
      \left(k^{\ast}(x_i,x_j)\right)_{i,j = 1}^n&
      \left(L_2(x_i,x_j)\right)_{i,j = 1}^n
    \end{pmatrix} : M_2(M_n(\mathcal{A})) \to M_2(M_n(\mathcal{B}))
  \end{equation*}
  is completely positive.
\item[(\emph{iv}\kern1pt)] Given any finite set $F = \{x_1,x_2,\ldots
  x_n\}$ of $X$ there exist
  $L_1,L_2\in{\mathbb{K}^{+}_F(\mathcal{A,B})}$ such that the map
  \begin{equation*}
    \begin{pmatrix}
      L_1(x_i,x_j) & k(x_i,x_j)\\
      k^{\ast}(x_i,x_j)& L_2(x_i,x_j)
    \end{pmatrix}_{i,j = 1}^n :
    M_n(M_2(\mathcal{A})) \to M_n(M_2(\mathcal{B}))
  \end{equation*}
  is completely positive.
\item[(\emph{v}\kern1pt)] Given any finite set $F = \{x_1,x_2,\ldots
  x_n\}$ of $X$ there exist
  $L_1,L_2\in{\mathbb{K}^{+}_F(\mathcal{A,B})}$ such that the map
  \begin{equation*}
    \begin{pmatrix}
      \left(L_1(x_i,x_j)\right)_{i,j = 1}^n &
      \left(k(x_i,x_j)\right)_{i,j = 1}^n \\[6pt]
      \left(k^{\ast}(x_i,x_j)\right)_{i,j = 1}^n &
      \left(L_2(x_i,x_j)\right)_{i,j = 1}^n
    \end{pmatrix}
    : M_2(M_n(\mathcal{A})) \to M_2(M_n(\mathcal{B}))
  \end{equation*}
  is completely positive.
\item[(\emph{vi}\kern1pt)] Given any finite set $F = \{x_1,x_2,\ldots
  x_n\}$ of $X$ there exist $P_1,P_2$ completely positive maps from
  $M_n(\mathcal{A})$ to $M_n(\mathcal{B})$ such that the map
  \begin{equation*}
    \begin{pmatrix}
      P_1 & S_{k_F}\\
      S_{k_F}^{\ast} & P_2
    \end{pmatrix}
    : M_2(M_n(\mathcal{A})) \to M_2(M_n(\mathcal{B}))
  \end{equation*}
  is completely positive, where
  \begin{equation*}
    S_{k_F} : M_n(\mathcal{A}) \to M_n(\mathcal{B}) :
    \left(a_{i,j}\right)_{i,j = 1}^n \mapsto
    \left(k(x_i,x_j)[a_{i,j}]\right)_{i,j = 1}^n
  \end{equation*}
  is the Schur product operator associated to the matrix
  $(k(x_i,x_j))$.
\end{itemize}

We prove under appropriate conditions on $\mathcal{B}$ and $k$ that
all six statements are in fact equivalent.  We begin by proving the
equivalence of the three ``global'' statements
(\emph{i}\kern1pt)--(\emph{iii}\kern1pt), the equivalence of the three
``local'' statements (\emph{iv}\kern1pt)--(\emph{vi}\kern1pt), and
that the global statements imply the local statements.  The proof
draws upon the following two results, the first of which is a routine
generalisation of the off-diagonal technique in~\cite[Theorem
8.3]{PA02}.

\begin{theorem}
  \label{thm:off-diagonal}
  Let the $C^*$-algebra $\mathcal{B}$ be injective.  Let $\mathcal{C}$
  be a unital $C^{\ast}$-subalgebra of both $\mathcal{A}$ and
  $\mathcal{B}$.  Let $\varphi:\mathcal{A}\to\mathcal{B}$ be a
  completely bounded, $\mathcal{C}$-bimodule map.  Then there exist
  completely positive $\mathcal{C}$-bimodule maps $\varphi_1,\varphi_2
  : \mathcal{A}\to \mathcal{B}$ with
  $\|\varphi_i\|_{cb}=\|\varphi\|_{cb}$ such that the map
  \begin{equation*}
    \Phi : M_2(\mathcal{A}) \to M_2(\mathcal{B}) :
    \begin{pmatrix}
      a & b \\ c & d
    \end{pmatrix}
    \mapsto
    \begin{pmatrix}
      \varphi_1(a) & \varphi(b) \\ \varphi^{\ast}(c) & \varphi_2(d)
    \end{pmatrix}
  \end{equation*}
  is completely positive.
\end{theorem}

\begin{lemma}
  \label{lem:bimodule_entrywise}
  Let $\varphi\in\mathcal{L}(M_n(\mathcal{A}),M_n(\mathcal{L(B)}))$.
  Then the following are equivalent.
  \begin{enumerate}
  \item $\varphi$ is a $\mathcal{D}_n$-bimodule map.
  \item For all $i,j=1,\ldots,n$ and all $A\in M_n(\mathcal{A})$ we
    have
    \begin{equation*}
      E_{i,j}\ast\varphi(A)=\varphi(E_{i,j}\ast{A}),
    \end{equation*}
    where $E_{i,j}$ is a matrix unit (that is, the $M_n$ element with
    $1$ in the $(i,j)^{\mathrm{th}}$ position and $0$ elsewhere) and
    $\ast$ is the entrywise (i.e, Schur) product.
  \item $\varphi$ acts entrywise on $M_n(\mathcal{A})$.
  \end{enumerate}
\end{lemma}

\begin{proof}
  It is clear that $(\emph{ii}\kern1pt)$ and $(\emph{iii}\kern1pt)$ of
  the lemma statement are equivalent.  The equivalence of
  $(\emph{i}\kern1pt)$ and $(\emph{ii}\kern1pt)$ follows from
  \begin{equation*}
    E_{i,i}AE_{j,j}=E_{i,j}\ast{A},
  \end{equation*}
  an easily checked equality.
\end{proof}

\begin{theorem}
  \label{thm:statement_equivalences}
  For the above statements, the following implications hold:
  \begin{equation*}
    \mathrm{(\text{i})} \Leftrightarrow \mathrm{(\text{ii})}
    \Leftrightarrow \mathrm{(\text{iii})} \Rightarrow
    \mathrm{(\text{iv})} \Leftrightarrow \mathrm{(\text{v})}
    \Leftrightarrow \mathrm{(\text{vi})}.
  \end{equation*}
\end{theorem}

\begin{proof}
  Statement (\emph{ii}\kern1pt) is a restatement of (\emph{i}\kern1pt)
  using a characterisation of completely positive kernels from
  \cite[Lemma 3.2.1]{BBLS03}.  The equivalence of (\emph{ii}\kern1pt)
  and (\emph{iii}\kern1pt) follows by using the (complete positivity
  preserving) canonical shuffle of matrices from \cite[Chapter
  8]{PA02}, and likewise for (\emph{iv}\kern1pt) and
  (\emph{v}\kern1pt).  Statement(\emph{iii}\kern1pt) implies
  (\emph{v}\kern1pt) by a restriction of kernels to a finite subset of
  $X$, and statement (\emph{v}\kern1pt) gives an explicit form for the
  completely positive maps $P_1,P_2$ in (\emph{vi}\kern1pt).  It is
  interesting to note that at this point we have not used the
  assumptions that $\mathcal{A}$ is unital and $\mathcal{B}$ is
  injective.

  Finally we prove (\emph{vi}\kern1pt) implies (\emph{v}\kern1pt).  By
  assumption there exist completely positive maps $P_1,P_2$ such that
  \begin{equation*}
    \begin{pmatrix}
      P_1 & S_{k_F} \\ S_{k_F}^{\ast} & P_2
    \end{pmatrix}
    \geq 0,
  \end{equation*}
  so the matrix is a completely bounded map.  This implies that
  $S_{k_F}$ is a completely bounded map.  It is a Schur product map,
  so acts entrywise.  By Lemma~\ref{lem:bimodule_entrywise} it is a
  $\mathcal{D}_n$-bimodule map.  Apply Theorem~\ref{thm:off-diagonal}
  and require that the completing maps are $\mathcal{D}_n$-bimodule
  maps, $P_1',P_2'$.  These are entrywise maps, so can be identified
  with elements of $M_n(\mathcal{L(A,B)})$, which we denote by
  $P_1'',P_2''$.  Define completely positive kernels $L_1,L_2$ on $F$
  by
  \begin{equation*}
    L_1(x_i,x_j)[a] := (P_1'')_{i,j}[a], \qquad L_2(x_i,x_j)[a] :=
    (P_2'')_{i,j}[a], \qquad i,j=1,\ldots, n.
  \end{equation*}
  These satisfy the conditions of statement (\emph{v}\kern1pt).
\end{proof}

\section{Topologies on kernel spaces}
\label{sec:topologies}

In this section, we consider topologies on the space
$\mathbb{K}_F(\mathcal A , \mathcal B)$ of all kernels on $F$ where
$F$ is a finite subset of $X$.  The topologies constructed will be
used in the last section to prove, under appropriate restrictions on
the space $\mathcal B$, the existence of Kolmogorov decompositions of
uniformly completely bounded kernels.

\subsection{The topology of bounded linear maps into $\mathcal{L(H)}$}
\label{subsec:topol-bound-linear-maps-into-LH}

Denote by $\mathcal{L}^1(\mathcal{H})$ the ideal of trace class
operators on $\mathcal{H}$.  Use the fact that $\mathcal{L(H)}$ is the
dual of $\mathcal{L}^1(\mathcal{H})$ to identify $\mathcal{L(A,L(H))}$
as an isometrically isomorphic copy of $(\mathcal{A}\otimes
\mathcal{L}^1(\mathcal{H}))^{\ast}$.  Given a $\varphi \in
\mathcal{L(A,L(H))}$ we associate with it a linear functional
$L_\varphi$ defined on elementary tensors by
\begin{equation*}
  L_\varphi(a\otimes{R})=\varphi(a)(R)=\mathrm{tr}\,(\varphi(a)R).
\end{equation*}
To a linear functional $L \in (\mathcal{A} \otimes
\mathcal{L}^1(\mathcal{H}))^{\ast}$ and an element $a \in \mathcal{A}$
we associate a linear functional $L^a$ on $
\mathcal{L}^1(\mathcal{H})$ by $L^a(R) = L(a\otimes{R})$.  Thus we can
define a bounded, linear map $\varphi_L : \mathcal{A} \to
\mathcal{L(H)}$ by $\varphi_L(a) = L^a$.  Clearly,
\begin{equation*}
  \varphi_{L_\varphi} = \varphi \quad\text{and}\quad L_{\varphi_L} = L.
\end{equation*}
The space $(\mathcal{A} \otimes \mathcal{L}^1(\mathcal{H}))^{\ast}$
carries a natural weak-$*$ topology.  Formally, we endow
$\mathcal{L(A,L(H))}$ with the weakest topology with respect to which
the isometric isomorphism $\varphi\mapsto L_\varphi$ is (weak-$\ast$)
continuous.  This is called the \textbf{bounded weak}, or \textbf{BW
  topology}, see \cite[Chapter 7]{PA02}.

We define a weaker topology, called the \textbf{bounded-bounded weak},
or \textbf{BBW topology}, to be the weakest topology with respect to
which the evaluations $E_{a,R}:\mathcal{L(A,L(H))}\to \mathbb{C}$
defined by $E_{a,R}(\varphi) = L_\varphi(a\otimes{R})$ are continuous
for all $a\in\mathcal{A}$ and $R \in \mathcal{L}^1(\mathcal{H})$.

\begin{proposition}
  \label{prop:BBW_is_pointwise_weak_star}
  A net $\varphi_\alpha$ in $\mathcal{L(A,L(H))}$ converges to
  $\varphi$ in the BBW topology if and only if $\varphi_\alpha(a)$
  converges weak-$\ast$ to $\varphi(a)$ for all $a\in\mathcal{A}$.
\end{proposition}

\begin{proof}
  Let $\varphi_\alpha$ be a net in $\mathcal{L(A,L(H))}$.  Then
  $\varphi_\alpha$ converges BBW to $\varphi$ if and only if
  $E_{a,R}\varphi_\alpha$ converges to $E_{a,R}\varphi$ for all
  $a\in\mathcal{A}$ and $R\in\mathcal{L}^1(\mathcal{H})$.  This
  happens if and only if $L_{\varphi_\alpha}(a\otimes{R})$ converges
  to $L_\varphi(a\otimes{R})$ which occurs if and only if
  $\varphi_\alpha(a)(R)$ converges to $\varphi(a)(R)$.  The last
  statement is equivalent to $\varphi_\alpha(a)$ converging
  weak-$\ast$ to $\varphi(a)$ for all $a\in\mathcal{A}$.
\end{proof}

We now restate some results from \cite{PA02} to affirm that useful
statements about the bounded weak topology remain true about the
bounded-bounded weak topology.

\begin{corollary}
  \label{cor:BBW_is_bounded_BW}
  A bounded net converges BBW if and only if it converges BW.
\end{corollary}

\begin{proof}
  By \cite[Lemma 7.2]{PA02}, if $\varphi_\alpha$ is a bounded net,
  $\varphi_\alpha(a)$ converges weak-$\ast$ to $\varphi(a)$ for all
  $a\in\mathcal{A}$ if and only if $\varphi_\alpha$ converges BW to
  $\varphi$.
\end{proof}

\begin{corollary}
  \label{cor:BBW_7pt3}
  A net $\varphi_\alpha$ in $\mathcal{L(A,L(H))}$ converges to
  $\varphi$ in the BBW topology if and only if
  $\ip{\varphi_\alpha(a)h}{k}$ converges to $\ip{\varphi(a)h}{k}$ for
  all $h,k\in\mathcal{H}$ and $a\in\mathcal{A}$.
\end{corollary}

\begin{proof}
 Combine the previous corollary with \cite[Proposition 7.3]{PA02}.
\end{proof}

\begin{proposition}
  \label{prop:bounded_BW_closed_is_BBW_closed}
  Any bounded, BW-closed subset $V$ of $\mathcal{L(A,L(H))}$ is
  BBW-closed.
\end{proposition}

\begin{proof}
  Let $v$ be in the BBW-closure of $V$.  Then there is a bounded
  BBW-convergent net, so a bounded BW-convergent net, converging to
  $v$.  Hence $v$ is in the BW-closure of the BW-closed set $V$, i.e,
  $v\in{V}$.
\end{proof}

\begin{proposition}
  \label{prop:bounded_BW_compact_is_BBW_compact}
  Any bounded, BW-compact subset $K$ of $\mathcal{L(A,L(H))}$ is
  BBW-compact.
\end{proposition}

\begin{proof}
  Let $\left\{\mathcal{U}_\lambda\right\}$
  be a BBW-open cover of $K$.  Then
  $\left\{\mathcal{U}_\lambda\right\}$ is a BW-open cover of $K$, so
  has a finite subcover, comprising BBW-open sets.
\end{proof}

\begin{proposition}
  \label{prop:BBW-top-is-Hausdorff}
  The space $\mathcal{L(A,L(H))}$ is Hausdorff in the BBW-topology.
\end{proposition}

\begin{proof}
  By virtue of being identified with the continuous linear functionals
  on $\mathcal{L(A,L(H))}$, the elements of $\mathcal{A} \otimes
  \mathcal{L}^1(\mathcal{H})$ separate points of
  $\mathcal{L(A,L(H))}$.  So let $\varphi_1\neq \varphi_2 \in
  \mathcal{L(A,L(H))}$, and suppose that $X \in \mathcal{A} \otimes
  \mathcal{L}^1(\mathcal{H})$ is such that $L_{\varphi_1} (X) \neq
  L_{\varphi_2} (X)$.  Finite linear combinations of elementary tensor
  products are norm dense in $\mathcal{A} \otimes
  \mathcal{L}^1(\mathcal{H})$, so by continuity, we may assume without
  loss of generality that $X = \sum_j^n a_j\otimes R_j$, where $a_j\in
  \mathcal A$ and $R_j \in \mathcal{L}^1(\mathcal{H})$ for all $j$.
  From this it follows that for some elementary tensor $a\otimes R$,
  $L_{\varphi_1} (a\otimes R) \neq L_{\varphi_2} (a\otimes R)$; that
  is, $E_{a,R}(\varphi_1) \neq E_{a,R}(\varphi_2)$.

  Set $\epsilon = \left|E_{a,R}(\varphi_1) -
    E_{a,R}(\varphi_2)\right|$.  For $j=1,2$, let $B_j$ be the ball of
  radius $\epsilon/3$ in $\mathbb C$ centred at
  $E_{a,R}(\varphi_j)$.  Observe that that the set $E_{a,R}^{-1}(B_j)$
  is open by continuity, contains $\varphi_j$ and satisfies
  $E_{a,R}^{-1}(B_1) \cap E_{a,R}^{-1}(B_2) =
  \emptyset$, finishing the proof.
\end{proof}

Let $\mathcal{D}_n$ be the algebra of diagonal, scalar-valued
matrices and let
$\mathcal{E}(M_n(\mathcal{A}),M_n(\mathcal{L(H)}))$ be the set of
$\mathcal{D}_n$-bimodule maps.  We shall need the following result
for later use.

\begin{theorem}
  \label{thm:bimodule_closed}
  The set $\mathcal{E}(M_n(\mathcal{A}),M_n(\mathcal{L(H)}))$ is BBW
  closed in $\mathcal{L}(M_n(\mathcal{A}),M_n(\mathcal{L(H)}))$.
\end{theorem}

\begin{proof}
  Let $\varphi_\alpha$ be a net in
  $\mathcal{E}(M_n(\mathcal{A}),M_n(\mathcal{L(H)}))$ converging to
  $\varphi\in\mathcal{L}(M_n(\mathcal{A}),M_n(\mathcal{L(H)}))$.  We
  define two maps
  \begin{equation*}
    \psi^{i,j}_1 : \mathcal{L}(M_n(\mathcal{A}),M_n(\mathcal{L(H)}))
    \to \mathcal{L}(M_n(\mathcal{A}),M_n(\mathcal{L(H)})) :
    \theta \mapsto E_{i,j}\ast\theta,
  \end{equation*}
  where $(E_{i,j}\ast\theta)(A):=E_{i,j}\ast\theta(A)$, and
  \begin{equation*}
    \psi^{i,j}_2 : \mathcal{L}(M_n(\mathcal{A}),M_n(\mathcal{L(H)}))
    \to \mathcal{L}(M_n(\mathcal{A}),M_n(\mathcal{L(H)})) :
    \theta \mapsto \theta\ast E_{i,j},
  \end{equation*}
  where $(\theta\ast E_{i,j})(A):=\theta(E_{i,j}\ast A)$.  By
  Lemma~\ref{lem:bimodule_entrywise}, elements of
  $\mathcal{E}(M_n(\mathcal{A}),M_n(\mathcal{L(H)}))$ are
  characterised by the property
  $\psi^{i,j}_1(\theta)=\psi^{i,j}_2(\theta)$ for all
  $i,j=1,\ldots,n$.  Suppose that each of these functions is
  BBW-continuous.  Then, since $\varphi_\alpha
  \stackrel{BBW}{\longrightarrow} \varphi$, it follows that
  $\psi^{i,j}_1(\varphi_\alpha) \stackrel{BBW}{\longrightarrow}
  \psi^{i,j}_1(\varphi)$ and $\psi^{i,j}_2(\varphi_\alpha)
  \stackrel{BBW}{\longrightarrow} \psi^{i,j}_2(\varphi)$ for all
  $i,j=1,\ldots,n$.  The convergent nets
  $\psi^{i,j}_1(\varphi_\alpha)$ and $\psi^{i,j}_2(\varphi_\alpha)$
  are identical, from which it follows that
  \begin{equation*}
    \psi^{i,j}_1(\varphi) = \psi^{i,j}_2(\varphi), \qquad
    i,j = 1, \ldots ,n.
  \end{equation*}
  That is, $\varphi \in
  \mathcal{E}(M_n(\mathcal{A}),M_n(\mathcal{L(H)}))$.

  It remains to prove that each $\psi^{i,j}_1,\psi^{i,j}_2$ is
  BBW-continuous.  Let us begin by observing that, by virtue of its
  definition, the weak-$*$ topology on $(M_n(\mathcal{A})\otimes
  M_n(\mathcal{L}^1(\mathcal{H})))^{\ast}$ is generated by basic open
  sets of the form
  \begin{equation*}
    \mathcal{U}'=\widehat{(A\otimes R)}^{-1}
    (\mathcal{B}_\epsilon(z_0)),
  \end{equation*}
  with $\mathcal{B}_\epsilon(z_0)$ a ball in $\mathbb{C}$ centred at
  $z_0$ of radius $\epsilon>0$, $A\in M_n(\mathcal{A})$, $R \in
  \mathcal{L}^1(\mathcal{H}^n) \cong M_n(\mathcal{L}^1(\mathcal{H}))$
  and $\widehat{(A\otimes R)}$ the evaluation function.  It follows
  that the BBW topology has basic open sets of the form
  \begin{equation*}
    \mathcal{U} =
    \left\{\varphi\in\mathcal{L}(M_n(\mathcal{A}),M_n(\mathcal{L(H)}))
      : \left|L_\varphi(A\otimes R)-z_0\right| < \epsilon \right\}.
  \end{equation*}
  Since the weak-$*$ continuous linear functionals separate
  $(M_n(\mathcal{A})\otimes M_n(\mathcal{L}^1(\mathcal{H})))^{\ast}$,
  there exists $\varphi_0$ such that $L_{\varphi_0}(A\otimes R) \neq
  0$.  Replacing $\varphi_0$ by $\tfrac{z_0}{L_{\varphi_0}(A\otimes
    R)} L_{\varphi_0}(A\otimes R)$, we then see that a basis for the
  BBW topology is generated by open sets
  \begin{equation*}
    \mathcal{U} = \left\{\varphi \in \mathcal{L}(M_n(\mathcal{A}),
      M_n(\mathcal{L(H)})) : \left|(L_\varphi - L_{\varphi_0})(A
        \otimes R)\right| < \epsilon \right\}
  \end{equation*}
  as $\varphi_0,A,R$ and $\epsilon$ vary.

  We now consider ${\psi^{i,j}_1}^{-1}(\mathcal{U})$, for each of
  these basic open sets.  Of course if $\mathcal{U}$ does not
  intersect the range of $\psi^{i,j}_1$ then
  ${\psi^{i,j}_1}^{-1}(\mathcal{U})=\emptyset$ is trivially open.  As
  above we consider open sets $\mathcal{U}$ centred at
  $L_{{\psi^{i,j}_1}(\varphi_0)}(A\otimes R)$ for some $\varphi_0$,
  $A$ and $R$.  Then
  \begin{equation*}
    {\psi^{i,j}_1}^{-1} \left(\left\{ \varphi \in
        \mathcal{L}(M_n(\mathcal{A}),M_n(\mathcal{L(H)})) :
        \left|\left(L_\varphi - L_{{\psi^{i,j}_1}(\varphi_0)}\right)
          (A \otimes R) \right| < \epsilon \right\}\right)
  \end{equation*}
  is the set
  \begin{equation*}
    \left\{ \varphi \in \mathcal{L}(M_n(\mathcal{A}),
      M_n(\mathcal{L(H)})) : \left| \left(L_{{\psi^{i,j}_1}(\varphi)}
          - L_{{\psi^{i,j}_1}(\varphi_0)}\right) (A\otimes R) \right|
      < \epsilon \right\}.
  \end{equation*}
  Now setting ${\psi^{i,j}_1}(\theta) = E_{i,j}\ast\theta$ and
  $L_\varphi(A \otimes R) = \mathrm{tr}(\varphi(A)R)$ we get that
  \begin{equation*}
    {\psi^{i,j}_1}^{-1}(\mathcal{U}) = \left\{ \varphi \in
      \mathcal{L}(M_n(\mathcal{A}), M_n(\mathcal{L(H)})) : \left|
        \mathrm{tr}\left([E_{i,j}\ast(\varphi - \varphi_0)(A)]R\right)
      \right| < \epsilon \right\}.
  \end{equation*}
  Making the identifications $M_n(\mathcal{L(H)}) =
  \mathcal{L}(\mathcal{H}^n)$ and $M_n(\mathcal{L}^1(\mathcal{H})) =
  \mathcal{L}^1(\mathcal{H}^n)$,
  \begin{equation*}
    [E_{i,j}\ast(\varphi - \varphi_0)(A)]R = [(\varphi - \varphi_0)(A)
    \ast E_{i,j}]R = [(\varphi - \varphi_0)(A)] [E_{i,j} \ast {R}] =
    [(\varphi - \varphi_0)(A)] \widetilde{R}
  \end{equation*}
  where $\widetilde{R} \in \mathcal{L}^1(\mathcal{H}^n)$.  Thus
  \begin{equation*}
    {\psi^{i,j}_1}^{-1}(\mathcal{U}) = \left\{ \varphi \in
      \mathcal{L}(M_n(\mathcal{A}), M_n(\mathcal{L(H)})) : \left|
        \mathrm{tr}\left((\varphi - \varphi_0)(A) \widetilde{R}\right)
      \right| < \epsilon \right\}
  \end{equation*}
  is clearly BBW-open.

  A substantially identical argument shows that, for $\mathcal{U}$ in
  the analogously chosen basis,
  \begin{equation*}
    {\psi^{i,j}_2}^{-1}(\mathcal{U}) = \left\{ \varphi \in
      \mathcal{L}(M_n(\mathcal{A}), M_n(\mathcal{L(H)})) : \left|
        \mathrm{tr}\left((\varphi - \varphi_0)(E_{i,j} \ast A)R\right)
      \right| < \epsilon\right\}.
  \end{equation*}
  Defining $\widetilde{A}:=E_{i,j}\ast A$, it becomes clear that
  ${\psi^{i,j}_2}^{-1}(\mathcal{U})$ is open.  Thus $\psi^{i,j}_1$ and
  $\psi^{i,j}_2$ are continuous, and the result is proven.
\end{proof}

Finally, some brief observations on relative topologies.

\begin{proposition}
  \label{prop:relative_topology_compactness}
  Let $\mathbb{X}$ be a space, $\mathbb{Y}$ a topological space with
  topology $\tau_\mathbb{Y}$, and let
  $f:\mathbb{X}\to\mathbb{Y}$ be injective.  Endow
  $\mathbb{X}$ with the weakest topology $\tau_\mathbb{X}$ such that
  $f$ is continuous.  Suppose that $A\subset\mathbb{X}$, so
  $f(A)\subset\mathbb{Y}$ and denote the relative topologies of $A$
  and $f(A)$
  \begin{equation*}
    \tau_A := \left\{ \mathcal{U} \cap {A} : \mathcal{U} \in
      \tau_\mathbb{X} \right\}, \qquad
    \tau_B := \left\{ \mathcal{U} \cap {f(A)} : \mathcal{U} \in
          \tau_\mathbb{Y} \right\}.
  \end{equation*}
  Then $\tau_A$ is the weakest topology with respect to which $f|_A$
  is ($\tau_B$-)continuous.  Further, if $f(A)$ is $\tau_B$-compact,
  then $A$ is $\tau_A$-compact.
\end{proposition}

\begin{proof}
  Let the weakest topology with respect to which $f|_A$ is
  ($\tau_B$-)continuous be $\tau_A'$.  Then
  \begin{equation*}
    \tau_A' := \left\{ f^{-1}(\mathcal{U} \cap {f(A)}) : \mathcal{U}
      \in \tau_\mathbb{Y} \right\} = \left\{ f^{-1}(\mathcal{U}) \cap
      {A} : \mathcal{U} \in \tau_\mathbb{Y} \right\} = \tau_A.
  \end{equation*}
  For compactness, simply note that $f|_A$ is a homeomorphism of
  $(A,\tau_A)$ and $(f(A),\tau_B)$.
\end{proof}

\begin{lemma}[{\cite[\S12, Corollary~1, p.~68]{HO75}}]
  Let $X$ be a normed, linear space.  Then every bounded subset of
  $X^{\ast}$ is relatively weak-$\ast$ compact.
\end{lemma}

\begin{corollary}
  \label{cor:bounded_sets_relatively_compact}
  Every bounded subset of $\mathcal{L(A,L(H))}$ is relatively
  BBW compact.
\end{corollary}

\begin{proof}
  By Proposition~\ref{prop:relative_topology_compactness}, the bounded
  subsets of $\mathcal{L(A,L(H))}$ are relatively BW-compact.  Let $A$
  be such a subset, and take a cover $\left\{\mathcal{U}_\alpha \cap
    {A}\right\}_\alpha$ for some collection
  $\left\{\mathcal{U}_\alpha\right\}_\alpha$ of BBW-open sets.  Since
  BBW-open sets are BW-open, the sets $\mathcal{U}_\alpha \cap A$ are
  relatively BW-open.  Relative BW-compact\-ness gives the finite
  subcover we require.
\end{proof}

\subsection{The pointwise $\sigma$-BBW topology on the
  kernels}
\label{subsec:pointw-sigma-bound-wk-top}

For $x,y \in X$, the evaluation maps $\mathbb{F}_{x,y} :
\mathbb{K}_X(\mathcal{A,B}) \to \mathcal{L(A,B)}$ are defined as
$\mathbb{F}_{x,y}(k) = k(x,y)$.  Let
$F=\left\{x_1,x_2,\ldots,x_n\right\}$, $x_i\in{X}$.  We define the
\textbf{pointwise $\sigma$-BBW topology} $\tau^p_F$ on
$\mathbb{K}_F(\mathcal{A,B})$ to be the weakest topology such that
for all $x,y\in{F}$ the maps
\begin{equation*}
  \mathbb{F}_{x,y} : \mathbb{K}_F(\mathcal{A,B}) \to
  \mathcal{L(A,B)}
\end{equation*}
are continuous, where $\mathcal{L(A,B)}$ is endowed with the BBW
topology.  It is then clear that when $G \subset F$, the restriction
maps
\begin{equation}
  \label{eq:1}
  \mathbb{K}_F(\mathcal{A,B}) \to\mathbb{K}_G(\mathcal{A,B})
  : k \mapsto k|_G
\end{equation}
are automatically continuous.  The evaluations $\mathbb{F}_{x,y}$
separate $\mathbb{K}_F(\mathcal{A,B})$ and the BBW topology is locally
convex and Hausdorff, so each $\tau^p_F$ is locally convex and
Hausdorff.

\subsection{The local $\sigma$-BBW topology on the kernels}
\label{subsec:local-sigma-bounded-wk-tops}

There is another topology worth considering on kernel spaces.  From
the previous discussion, for a fixed
$F=\left\{x_1,x_2,\ldots,x_n\right\}$ there is an identification
\begin{equation*}
  j_F:\mathbb{K}_F(\mathcal{A,B})\to
  \mathcal{L}(M_n(\mathcal{A}),\mathcal{L}(\mathcal{H}^n)): k\mapsto
  \widehat{\sigma}\left(S_{\left(k(x_i,x_j)\right)}\right)
\end{equation*}
of a kernel with the Schur product operator associated to the matrix
$\left(k(x_i,x_j)\right),x_i\in{F}$.  We define the \textbf{local
  $\sigma$-BBW topology} $\tau^\ell_F$ on
$\mathbb{K}_F(\mathcal{A,B})$ to be the weakest topology such that
$j_F$ is continuous.  The map $j_F$ is injective and the BBW topology
is a locally convex Hausdorff topology.  Hence the local $\sigma$-BBW
topology is also locally convex and Hausdorff.  In general we will
abuse notation regarding Schur product operators, writing either
$(k(x_i,x_j))_F$ or $S_{k_F}$ for the map
$S_{\left(k(x_i,x_j)\right)}$.

\subsection{Equivalence of the two $\sigma$-BBW topologies}
\label{subsec:equiv-two-topol}

\begin{lemma}
  \label{lem:equivalence_of_topologies}
  For a given faithful unital $*$-representation $\sigma$ of the
  unital $C^*$-algebra $\mathcal B$, the pointwise $\sigma$-BBW
  topology on $\mathbb{K}_F$ is the same as the local $\sigma$-BBW
  topology on $\mathbb{K}_F$.
\end{lemma}

\begin{proof}
  Since the topologies are defined as being the weakest making certain
  maps continuous, it suffices to show that these maps are continuous
  in both topologies, and this is done by showing that if a net
  $(k_\alpha)$ of kernels in $\mathbb{K}_F$ converges to $k$ in one of
  the topologies, it does so in the other.

  So assume that $(k_\alpha)$ is a net of kernels in $\mathbb{K}_F$
  converging to $k$ in the local $\sigma$-BBW topology $\tau^\ell_F$.
  Then for any $\tilde a \in M_n(\mathcal{A})$ and $\tilde R\in
  \mathcal{L}^1(\mathcal{H}^n)\cong M_n(\mathcal{L}^1(\mathcal{H}))$,
  \begin{equation*}
    \mathrm{tr}\left(\widehat\sigma(k_\alpha(x_i,x_j))[\tilde a]
      \tilde R\right) \to
    \mathrm{tr}\left(\widehat{\sigma}(k(x_i,x_j))[\tilde a]
      \tilde R\right).
  \end{equation*}
  In particular, if we fix $i,j$ and choose $\tilde a = a\otimes
  E_{ij}$ and $\tilde R = R\otimes E_{ji}$, where $a\in\mathcal A$,
  $R\in \mathcal{L}^1(\mathcal{H})$, and $E_{ij}, E_{ji}$ are matrix
  units in $M_n(\mathbb C)$, we find that
  \begin{equation*}
    \mathrm{tr}\left(\sigma((k_\alpha(x_i,x_j)[a])R\otimes E_{ii}
    \right) \to
    \mathrm{tr}\left(\sigma((k(x_i,x_j)[a])R\otimes E_{ii}
    \right).
  \end{equation*}
  Thus $\mathrm{tr}\left(\sigma(k_\alpha(x_i,x_j)[a])R\right) \to
  \mathrm{tr}\left(\sigma(k(x_i,x_j)[a])R\right),$ from which it
  follows that $\mathbb F_{x_i,x_j}$ is continuous.  The choice of $i$
  and $j$ are arbitrary, hence $\tau^p_F \subseteq \tau^\ell_F$.

  On the other hand, suppose that $(k_\alpha)$ converges to $k$ in the
  pointwise $\sigma$-BBW topology $\tau^p_F$.  Then for $\tilde R =
  (R_{ij}) \in M_n(\mathcal{L}^1(\mathcal{H}))$, $\tilde a = (a_{ij})
  \in M_n(\mathcal{A})$ and any $\alpha$,
  \begin{equation}
    \label{eq:2}
    \mathrm{tr}\left(\widehat{\sigma}(k_\alpha(x_i,x_j))[\tilde a]
      \tilde R\right)
  \end{equation}
  is just the sum of the traces of the diagonal elements of
  $\widehat{\sigma}(k_\alpha(x_i,x_j))[\tilde a]\tilde R$, which are
  all linear combinations of things of the form
  $\sigma(k_\alpha(x,y)[a])R$, which by definition of the topology
  $\tau^p_F$, converge to the same linear combinations, with
  $k_\alpha$ replaced by $k$.  Consequently, \eqref{eq:2} converges to
  $\mathrm{tr}\left(\widehat{\sigma}\left(k(x_i,x_j)[\tilde a]\right)
    \tilde R\right)$.  We conclude that $\tau^\ell_F \subseteq
  \tau^p_F$, and therefore the two topologies agree.
\end{proof}

Henceforth we denote these two equivalent topologies on $\mathbb{K}_F
(\mathcal A , \mathcal B)$ by $\tau_F$.  We have already observed that
the restriction maps in \eqref{eq:1} are continuous when
${\mathbb{K}}_F ( \mathcal A , \mathcal B)$ has the $\tau_F$ topology.
We endow a ball of completely bounded kernels
\begin{equation*}
  \mathbb{K}_F^r(\mathcal{A,B}) := \left\{k \in
    \mathbb{K}_F(\mathcal{A,B}) : \left\| (k(x_i,x_j))_F \right\|_{cb}
    \leq r \right\}
\end{equation*}
with a relative $\tau_F$-topology, denoted $\tau_F^r$.

\begin{lemma}
  \label{lem:balls-cpt}
  Let $\mathcal B$ be unital.  The balls
  \begin{equation*}
    \mathbb{K}_F^r(\mathcal{A,B}) :=
    \left\{k\in\mathbb{K}_F(\mathcal{A,B}):
      {\left\|\left(k(x_i,x_j)\right)_F\right\|}_{cb}
      \leq{r}\right\}
  \end{equation*}
  are $\tau_F^r$-compact.
\end{lemma}

\begin{proof}
  The set $j_F\left(\mathbb{K}_F^r(\mathcal{A,B})\right) =
  j_F\left(\mathbb{K}_F(\mathcal{A,B})\right) \cap \mathrm{CB}^r_F$ is
  bounded, so by Corollary~\ref{cor:bounded_sets_relatively_compact}
  it is relatively BBW-compact.  Endow $\mathbb{K}_F^r(\mathcal{A,B})$
  with the weakest topology such that the restricted map
  \begin{equation*}
    j_F^r : \mathbb{K}_F^r(\mathcal{A,B}) \to
    j_F(\mathbb{K}_F^r(\mathcal{A,B})) : k \mapsto j_F(k)
  \end{equation*}
  is continuous.  By
  Proposition~\ref{prop:relative_topology_compactness}, that topology
  is $\tau_F^r$, $j_F^r$ is a homeomorphism, and
  $\mathbb{K}_F^r(\mathcal{A,B})$ is compact in that topology.
\end{proof}

\begin{proposition}
  \label{prop:positive_cone_compact}
  The truncated positive cone
  \begin{equation*}
    \mathbb{K}_F^{r,+}(\mathcal{A,B}) :=
    \left\{k\in\mathbb{K}_F^{+}(\mathcal{A,B}) :
      \left\|(k(x_i,x_j))_F\right\|_{cb} \leq r \right\}
  \end{equation*}
  is $\tau_F^r$-compact.
\end{proposition}

\begin{proof}
  By \cite[Proposition 7.4]{PA02} the truncated cone
  $CP^r(M_n(\mathcal{A}),M_n(\mathcal{L(H)}))$ of positive elements in
  $\mathcal{L}(M_n(\mathcal{A}),M_n(\mathcal{L(H)}))$ is closed.  The
  map $j_F^r$ is continuous, so the pre-image of that cone, the
  truncated cone of kernels, is $\tau_F^r$-closed in the
  $\tau_F^r$-compact ball of completely bounded kernels, so is itself
  $\tau_F^r$-compact.
\end{proof}

Now we obtain a characterisation of convergence in $\tau_F^r$ in the
manner of ~\cite[Proposition 7.3]{PA02}, or
Corollary~\ref{cor:BBW_7pt3}.

\begin{lemma}
  \label{lem:relative_convergence}
  Let $\mathbb{X}$ be a topological space, and $A\subset\mathbb{X}$ a
  closed subset.  Then a net $\left\{a_\alpha\right\}\subset{A}$
  converges to $a\in{A}$ in the relative topology if and only if
  $a_\alpha\to{a}$ in $\mathbb{X}$.
\end{lemma}

\begin{proof}
  Suppose $a_\alpha\to{a}$ in the relative topology.  Take an open
  neighbourhood $\mathcal{U}$ of $a$, in $\mathbb{X}$, so
  $A\cap\mathcal{U}$ is a relatively open neighbourhood of $a$ in $A$.
  Then by our supposition $a_\alpha$ is eventually in
  $A\cap\mathcal{U}$, so is eventually in $\mathcal{U}$.  Therefore
  $a_\alpha\to{a}$ in $\mathbb{X}$.  Conversely, let $a_\alpha\to{a}$
  in $\mathbb{X}$.  Since $A$ is closed, $a\in{A}$.  Take a relatively
  open neighbourhood of $a$, $A\cap\mathcal{U}$, $\mathcal{U}$ open in
  $\mathbb{X}$.  Then $a_\alpha$ is always in $A$ and eventually in
  $\mathcal{U}$, so is eventually in $A\cap\mathcal{U}$.  Thus
  $a_\alpha\to{a}$ relatively.
\end{proof}

\begin{proposition}
  \label{prop:relative_convergence_two}
  Let $k_\alpha$ be a net in $\mathbb{K}_F^r(\mathcal{A,B})$.  Let
  $\sigma:\mathcal{B}\to\mathcal{L(H)}$ be a $\ast$-representation,
  and allow $\sigma$ to induce $j_F,j_F^r$ and $\tau_F^r$ as described
  previously.  Then $k_\alpha\to k$ in $\tau_F^r$ if and only if for
  all $\left(a_{i,j}\right)\in M_n(\mathcal{A})$,
  $\left\{h_i\right\}_{i=1}^n, \left\{k_i\right\}_{i=1}^n \subset
  \mathcal{H}$ we have
  \begin{equation*}
    \ip{\oplus_i k_i}{\sigma^{(n)}\left((k_\alpha(x_i,x_j))_F
        \left(a_{i,j}\right)\right) \left(\oplus_i h_i\right)} \to
        \ip{\oplus_i k_i}{\sigma^{(n)}\left((k(x_i,x_j))_F
        \left(a_{i,j}\right) \right) \left(\oplus_i h_i\right)}.
  \end{equation*}
\end{proposition}

\begin{proof}
  Since $j_F^r$ is a homeomorphism, $k_\alpha\stackrel{\tau_F^r}\to k$
  if and only if $j_F^r(k_\alpha)\to j_F^r(k)$ in the relative BBW
  topology.  By Lemma~\ref{lem:relative_convergence}, this is
  equivalent to $j_F(k_\alpha)\to j_F(k)$ in the BBW topology.  The
  proof then follows from Corollary~\ref{cor:BBW_7pt3} and the
  definition of $j_F$.
\end{proof}

\section{Decomposition of regular completely bounded kernels}
\label{sec:decomp-reg-cb-kers}

As we noted in the remarks following the proof of
Theorem~\ref{thm:equiv_cond}, any kernel which has a Kolmogorov
decomposition must be a regular completely bounded kernel.  We now
draw on the topological results of the last section to prove that this
condition is both necessary and sufficient for $\mathcal{L(A,B)}$
valued kernels when $\mathcal{B}$ is injective.

\begin{theorem}
  \label{thm:reg-cb-kernels_have_K-decomps}
  If the $C^{\ast}$-algebra $\mathcal{B}$ is injective, then any
  regular completely bounded kernel $k\in\mathbb{K}_X(\mathcal{A,B})$
  has a Kolmogorov decomposition.
\end{theorem}

\begin{proof}
  Assume that $k\in\mathbb{K}_X(\mathcal{A,B})$ is a regular
  completely bounded kernel.  Then we can scale $k$ so that it is
  completely contractive, and so it suffices to prove the theorem
  under the assumption that $k$ is of this form.

  The plan for the proof is to invoke Lemma
  \ref{lem:kolm_decomp_characterisation}.  In other words, we shall
  show that the given kernel $k$ has a Kolmogorov decomposition by
  showing that there exist two completely positive kernels $L_1$ and
  $L_2$ in ${\mathbb K}_X(\cla , \clb)$ such that the kernel
 \begin{equation}
   \label{global-completion}
    \begin{pmatrix}
      L_1 & k \\ k^{\ast} & L_2
    \end{pmatrix} \in {\mathbb K}^+_X(\cla , M_2(\clb)).
  \end{equation}
  By complete boundedness of the given kernel $k$ and
  Theorem~\ref{thm:statement_equivalences}, this is known if $X$ is a
  finite set.  The global statement will be obtained by a form of
  transfinite induction and that is where the uniform bound will be
  needed.

  Well order the set $X$ and consider the directed set $\Lambda$ of
  finite strictly increasing sequences $\left\{\lambda_i\right\}$ of
  elements of $X$ with the order relation $\alpha\leq\beta$ in
  $\Lambda$ whenever $\alpha$ is a subsequence of $\beta$ (that is,
  there exists an injective, order-preserving map $ \widetilde{\alpha}
  : \left\{1,2, \ldots, |\alpha|\right\} \to \beta $ such that
  $\widetilde{\alpha}(j)=\alpha_j$).

  Recall that $ \mathbb{K}_F^1(\mathcal{A,B})$ is the
  $\tau_F^1$-compact set of all $k \in \mathbb{K}_F(\mathcal{A,B})$
  such that ${\| \left( k(x_i,x_j) \right)_F \|}_{cb} \leq{1}$, while
  $ \mathbb{K}_F^{1,+}(\mathcal{A,B}) = \mathbb{K}_F^1(\mathcal{A,B})
  \cap \mathbb{K}_F^+(\mathcal{A,B})$.  We define the space of local
  solutions by
  \begin{equation*}
    \mathbb{L}_F = \left\{ (L_1,L_2) : L_1, L_2 \in
    \mathbb{K}_F^{1,+}(\mathcal{A,B}) \text{ and }
       \begin{pmatrix}
        L_1 & k|_{F} \\ k^{\ast}|_{F} & L_2
      \end{pmatrix}
      \in \mathbb{K}_F^{1,+}(\mathcal{A},M_2(\mathcal{B})) \right\}.
  \end{equation*}
  By assumption this set is non-empty.  Thus we have a family of
  topological spaces $\{ {\mathbb L}_F : F \in \Lambda\}$ indexed by
  $\Lambda$.  For $G \geq F$ in $\Lambda$, let $f_{F,G} : \mathbb{L}_G
  \to \mathbb{L}_F$ be the restriction mapping; that is,
  $f_{F,G}(L_1,L_2)$ is the restriction of the kernels $L_1$ and $L_2$
  (a priori defined on $G$) to the smaller set $F$.  By definition of
  the topology $\tau_F^1$, the restriction maps $f_{F,G}$ are
  automatically continuous.  The topological spaces ${\mathbb L}_F$
  along with these maps $f_{F,G}$ comprise an \textsl{inverse system}
  and the existence of $L_1$ and $L_2$ satisfying
  \ref{global-completion} which we are aiming for is then equivalent
  to this inverse system having a non-empty \textsl{inverse limit}.
  This is guaranteed provided the so-called factor spaces ${\mathbb
    L}_F$ are compact and Hausdorff, which we prove below in detail.

  A collection of completely contractive, completely positive kernels
  supported on $F$ is, by Proposition
  \ref{prop:positive_cone_compact}, $\tau_F^1$-compact.  Consequently
  the product of two truncated positive cones,
  $\mathbb{K}_F^1(\mathcal{A,B})\times\mathbb{K}_F^1(\mathcal{A,B})$,
  is a compact Hausdorff in the product topology
  $\tau_F^1\times\tau_F^1$.  Thus, to show that $\mathbb{L}_F$ is
  $\tau_F^1\times\tau_F^1$-compact requires only that it be closed.

  Let $(L_1^\alpha,L_2^\alpha)$ be a convergent net in $\mathbb{L}_F$.
  By compactness of the positive cones, its limit is $(L_1,L_2)$ where
  $L_1,L_2$ are completely positive kernels, with
  $\|(L_1(x_i,x_j))\|_{cb},\|(L_2(x_i,x_j))\|_{cb}\leq{1}$.  By
  Proposition~\ref{prop:relative_convergence_two},
  $L^\alpha\stackrel{\tau_F^1}{\to}L$ if and only if for all
  $\left(a_{i,j}\right)\in M_n(\mathcal{A})$,
  $\left\{h_i\right\}_{i=1}^n, \left\{k_i\right\}_{i=1}^n \subset
  \mathcal{H}$ we have
  \begin{equation*}
    \ip{\oplus_i k_i}{\sigma^{(n)}\left((L_\alpha(x_i,x_j))_F
        \left(a_{i,j}\right)\right) \left(\oplus_i
        h_i\right)} \to
    \ip{\oplus_i k_i}{\sigma^{(n)}\left((L(x_i,x_j))_F
        \left(a_{i,j}\right)\right) \left(\oplus_i
        h_i\right)}.
  \end{equation*}
  View the space of $2\times{2}$ matrices with entries from
  $\mathbb{K}_F^1(\mathcal{A,B})$ as sitting inside
  $\mathbb{K}_F^{4}(\mathcal{A},M_2(\mathcal{B}))$ and endow that set
  with a relative topology induced by the representation $\sigma^{(2)}
  : M_2(\mathcal{B}) \to \mathcal{L}(\mathcal{H} \oplus \mathcal{H})$.
  Call this topology $\tau$.  The set of positive elements is
  $\tau$-closed (cf.~Proposition~\ref{prop:positive_cone_compact}) and
  we have a characterisation of convergence
  (cf.~Proposition~\ref{prop:relative_convergence_two}).

  Since $(L_1^\alpha,L_2^\alpha)\to(L_1,L_2)$, for all choices of
  $(a_{i,j})\in M_n(\mathcal{A})$ and
  $(h_i'),(k_i'),(h_i''),(k_i'')\in\oplus_{i=1}^n\mathcal{H}$ as
  above, we have
  \begin{equation*}
    \begin{split}
      & \ip{\oplus_i k_i'}{\sigma^{(n)}((L_1^\alpha(x_i,x_j))
        [a_{i,j}])(\oplus_i h_i')}
      + \ip{\oplus_i k_i''}{\sigma^{(n)}((L_2^\alpha(x_i,x_j))
        [a_{i,j}])(\oplus_i h_i'')} \\[3pt]
      \to\; & \ip{\oplus_i k_i'}{\sigma^{(n)}((L_1(x_i,x_j))
        [a_{i,j}])(\oplus_i h_i')}
      + \ip{\oplus_i k_i''}{\sigma^{(n)}((L_2(x_i,x_j))
        [a_{i,j}])(\oplus_i h_i'')} \\
    \end{split}
  \end{equation*}
  Add fixed terms to both sides,
  \begin{equation*}
    \begin{split}
      & \ip{\oplus_i k_i'}{\sigma^{(n)}((L_1^\alpha(x_i,x_j))
        [a_{i,j}])(\oplus_i h_i')}
      + \ip{\oplus_i k_i''}{\sigma^{(n)}((L_2^\alpha(x_i,x_j))
        [a_{i,j}])(\oplus_i h_i'')} \\
      &\qquad + \ip{\oplus_i k_i'}{\sigma^{(n)}((k(x_i,x_j))
        [a_{i,j}](\oplus_i h_i'')}
        + \ip{\oplus_i k_i''}{\sigma^{(n)}((k^{\ast}(x_i,x_j))
          [a_{i,j}])(\oplus_i h_i')} \\[3pt]
      \to\; & \ip{\oplus_i k_i'}{\sigma^{(n)}((L_1(x_i,x_j))
        [a_{i,j}])(\oplus_i h_i')}
      + \ip{\oplus_i k_i''}{\sigma^{(n)}((L_2(x_i,x_j))
        [a_{i,j}])(\oplus_i h_i'')} \\
      &\qquad + \ip{\oplus_i k_i'}{\sigma^{(n)}((k(x_i,x_j))
        [a_{i,j}](\oplus_i h_i'')}
        + \ip{\oplus_i k_i''}{\sigma^{(n)}((k^{\ast}(x_i,x_j))
          [a_{i,j}])(\oplus_i h_i')} \\
    \end{split}
  \end{equation*}
  and rewrite this as
  \begin{equation*}
    \begin{split}
      & \ip{
          \begin{pmatrix}
            \oplus_{i=1}^n k_i'\\[3pt] \oplus_{i=1}^n k_i''
          \end{pmatrix}
        }{
        \begin{pmatrix}
          \sigma^{(n)}((L_1^\alpha(x_i,x_j))[a_{i,j}]) &
          \sigma^{(n)}((k(x_i,x_j))[a_{i,j}]) \\[3pt]
          \sigma^{(n)}((k^{\ast}(x_i,x_j))[a_{i,j}]) &
          \sigma^{(n)}((L_2^\alpha(x_i,x_j))[a_{i,j}])
        \end{pmatrix}
        \begin{pmatrix}
          \oplus_{i=1}^n h_i'\\[3pt] \oplus_{i=1}^n h_i''
        \end{pmatrix}
        } \\[3pt]
      \to\; &
      \ip{
          \begin{pmatrix}
            \oplus_{i=1}^n k_i'\\[3pt] \oplus_{i=1}^n k_i''
          \end{pmatrix}
        }{
        \begin{pmatrix}
          \sigma^{(n)}((L_1(x_i,x_j))[a_{i,j}]) &
          \sigma^{(n)}((k(x_i,x_j))[a_{i,j}]) \\[3pt]
          \sigma^{(n)}((k^{\ast}(x_i,x_j))[a_{i,j}]) &
          \sigma^{(n)}((L_2(x_i,x_j))[a_{i,j}])
        \end{pmatrix}
        \begin{pmatrix}
          \oplus_{i=1}^n h_i'\\[3pt] \oplus_{i=1}^n h_i''
        \end{pmatrix}
        }. \\
    \end{split}
  \end{equation*}
  Perform a canonical shuffle on the large matrix and relabel the
  Hilbert space elements to conclude that for all choices of
  $(a_{i,j}),(h_i)$ and $(k_i)$,
  \begin{equation*}
    \begin{split}
      & \ip{\oplus_{i=1}^{2n} k_i}{\sigma^{(2n)}\left(\left(
            \begin{pmatrix}
              L_1^\alpha(x_i,x_j) & k(x_i,x_j) \\[3pt]
              k^{\ast}(x_i,x_j) & L_2^\alpha(x_i,x_j)
            \end{pmatrix}
          \right)
          \left[\left(a_{i,j}\right)\right]
        \right)
        \left(\oplus_{i=1}^{2n} h_i\right)}
      \\
      \to\; &
      \ip{\oplus_{i=1}^{2n} k_i}{\sigma^{(2n)}\left(\left(
            \begin{pmatrix}
              L_1(x_i,x_j) & k(x_i,x_j) \\[3pt]
              k^{\ast}(x_i,x_j) & L_2(x_i,x_j)
            \end{pmatrix}
          \right)
          \left[\left(a_{i,j}\right)\right]
        \right)
        \left(\oplus_{i=1}^{2n} h_i\right)}
      \\
    \end{split}
  \end{equation*}
  It follows from Proposition~\ref{prop:relative_convergence_two} that
  \begin{equation}
    \label{eq:3}
    \begin{pmatrix}
      L_1^\alpha & k|_F \\[3pt] k|_F^{\ast} & L_2^{\alpha}
    \end{pmatrix}
    \;\stackrel{\tau}{\to}\;
    \begin{pmatrix}
      L_1 & k|_F \\[3pt] k|_F^{\ast} & L_2
    \end{pmatrix}.
  \end{equation}

  The net on the left of \eqref{eq:3} is composed of completely
  positive kernels, and they are $\tau$-closed, so the limit is a
  completely positive kernel.  Consequently, $(L_1,L_2)$ belongs to
  the set of solutions $\mathbb{L}_F$, and so $\mathbb{L}_F$ is
  closed.  Thus it is a non-empty compact Hausdorff space for each
  $F$.

  We can now conclude by \cite[Theorem~6.B.11]{CV77} that the inverse
  limit system is non-empty; that is, there exists
  \begin{equation*}
    ((L_1^\lambda,L_2^\lambda)) \in \prod_{\lambda\in\Lambda}
    \mathbb{L}_\lambda \qquad\text{and} \qquad
    f_{\lambda_1,\lambda_2}((L_1^{\lambda_2}, L_2^{\lambda_2})) =
    (L_1^{\lambda_1},L_2^{\lambda_1}) \quad\text{whenever }\lambda_1
    \leq \lambda_2.
  \end{equation*}
  It follows that this object uniquely specifies an element of
  $\mathbb{K}_X^+(\mathcal{A,B})$ giving us a global completely
  positive matrix completion.  From
  Lemma~\ref{lem:kolm_decomp_characterisation} we deduce that the
  kernel $k$ therefore has a Kolmogorov decomposition.
\end{proof}

\section{Some examples}
\label{sec:some-examples}

We begin this section by giving an alternative characterisation of the
regularity condition.  This is used in proving that a particular
hermitian completely bounded kernel is \emph{not} regular, and hence
does not have a Kolmogorov decomposition.  We also consider some
examples where the regularity condition is automatic.

\begin{definition}
  \label{def:unif-regularity}
  We say that a kernel $k \in \mathbb{K}_X(\mathcal{A,B})$ is
  \textbf{uniformly regular} if it is regular and there exists
  $\alpha\in (0,1]$ and a choice of $c:X\to [\alpha,1]$ such that
  $\tilde{k}(x,y)[a] := c(x) k(x,y)[a] \overline{c(y)}$ is completely
  contractive.
\end{definition}

\begin{lemma}
  \label{lem:chain-unif-reg-implies-reg}
  Let $k \in \mathbb{K}_X(\mathcal{A,B})$, and let $\mathcal C$ be the
  collection of all subsets of $X$ with the property that if $Y\in
  \mathcal C$, then $k|Y$ is uniformly regular, partially ordered by
  inclusion.  There exists a chain $\widecheck{\mathcal C} \subset
  \mathcal C$ such that
  \begin{equation*}
    \bigcup_{Y\in \widecheck{\mathcal C}} Y = X
  \end{equation*}
  if and only if $k$ is regular.
\end{lemma}

\begin{proof}
  Suppose that there exists such a chain.  Then there exists $Y_0\in
  \widecheck{\mathcal C}$ such that $\mathrm{card}\,(X\backslash Y_0)
  \leq \omega_0$.  If there are only finitely many $Y$ in
  $\widecheck{\mathcal C}$ with $Y\geq Y_0$, then $X$ is a maximal
  element of $\widecheck{\mathcal C}$ and we are done.

  So we may assume that $\mathrm{card}\,(X\backslash Y_0) = \omega_0$
  and that there are countably infinitely many terms $Y_0 \leq Y_1
  \leq Y_2 \leq \ldots$ in $\widecheck{\mathcal C}$.  Since by
  assumption each $Y_k$ is uniformly regular, there is for each $k$ an
  $\alpha_{Y_k} > 0$ and $c_{Y_k}:Y_k \to [\alpha_{Y_k},1]$ so that
  the completely bounded norm of $(c_{Y_k}(x) k(x,y)[a]
  \overline{c_{Y_k}(y)})$ is at most $1$.  Without loss of generality,
  we may take $c_{Y_k}(x) = \alpha_{Y_k}$ for all $x\in Y_k$.

  We digress to observe that if for some $\epsilon, \beta \in (0,1)$,
  $\|A\| \leq 1-\epsilon$ and $\beta \left\| \begin{pmatrix} A & B \\
      C & D \end{pmatrix}\right\| \leq 1$, then
  \begin{equation*}
    \left\|
      \begin{pmatrix}
        A & \tfrac{\beta\epsilon}{4} B \\ \tfrac{\beta\epsilon}{4}
        C & (\tfrac{\beta\epsilon}{4})^2 D 
      \end{pmatrix}\right\| \\
    \leq
    (1-\tfrac{\beta\epsilon}{4}) \|A\| + \tfrac{\beta\epsilon}{4}
    \begin{pmatrix}
      A & B \\ C & D
    \end{pmatrix}
    \leq 1-\epsilon/2.
  \end{equation*}

  We use this observation to define $c:X\to (0,1]$ as follows.  For
  $x\in Y_0$, set $c(x) = \tfrac{1}{2}\alpha_{Y_0}$, so that we have
  $\| (c(x) k(x,y)[\cdot] \overline{c(y)})|_{Y_0}\|_{cb} \leq 1/2$.
  Next, for $x\in Y_1\backslash Y_0$, set $c(x) =
  \tfrac{1}{8}(\alpha_{Y_1})^2$, so that $\| (c(x) k(x,y)[\cdot]
  \overline{c(y)})|_{Y_1}\|_{cb} \leq 3/4$, and in general, for $x\in
  Y_{n+1}\backslash Y_n$, set $c(x) =
  \tfrac{1}{2^{2n+1}}(\alpha_{Y_{n+1}})^2$, so that $\| (c(x)
  k(x,y)[\cdot] \overline{c(y)})|_{Y_{n+1}}\|_{cb} \leq 1 - 2^{n+1}$.
  Since any finite set $F$ is contained in some $Y_n$, and $\| (c(x)
  k(x,y)[\cdot] \overline{c(y)})|_{F}\|_{cb}\| \leq (c(x)
  k(x,y)[\cdot] \overline{c(y)})|_{Y_{n}}\|_{cb} \leq 1$, we conclude
  that the kernel $k$ is regular.

  On the other hand, if $k$ is regular, then for $n\in \mathbb N$, we
  define
  \begin{equation*}
    Y_n := \{x\in X : c(x) \geq 1/n \}.
  \end{equation*}
  Then $\widecheck{\mathcal C} = \bigcup_n Y_n$ is a chain in
  $\mathcal C$ with $\bigcup_{Y\in \widecheck{\mathcal C}} Y = X$.
\end{proof}

\begin{corollary}
  \label{cor:cb-and-countable-ind-implies-decomp}
  If the $C^{\ast}$-algebra $\mathcal{B}$ is injective and the index
  set $X$ is countable, then any completely bounded kernel
  $k\in\mathbb{K}_X(\mathcal{A,B})$ has a Kolmogorov decomposition.
\end{corollary}

\begin{proof}
  Without loss of generality we take $X = \mathbb N$, and $Y_n =
  \{1,\ldots,n\}$, $n=1,2,\ldots$.  Then $k|_{Y_n}$ is uniformly
  regular, with $\alpha_{Y_k} = 1/\sqrt{\|k|_{Y_n}\|_{cb}}$ or $1$ if
  $\|k|_{Y_n}\|_{cb} < 1$.  Hence $\widecheck{\mathcal C} = \bigcup_n
  Y_n$ is a chain of sets such that each $k|_{Y_n}$ is uniformly
  regular and $\bigcup_n Y_n = X$, it follows from
  Lemma~\ref{lem:chain-unif-reg-implies-reg} that $k$ is regular.
  Hence by Theorem~\ref{thm:reg-cb-kernels_have_K-decomps} that $k$
  has a Kolmogorov decomposition.
\end{proof}

We now give an example, due to Scott McCullough, of a hermitian
completely bounded kernel which is not regular, and hence does not
have a Kolmogorov decomposition.

\begin{example}
  \label{exmpl:nonregular-cb-kernel}
  In this example we choose $X = [0,1]$ and $\mathcal A = \mathcal B =
  \mathbb C$.  In this case the entries of $k\in
  \mathbb{K}_X(\mathcal{A,B})$ are scalars and bounded implies
  completely bounded.  Define a completely bounded kernel by
  \begin{equation*}
    k(x,y) =
    \begin{cases}
      \dfrac{1}{|y-x|} & y\neq x \\[10pt]
      0 & y=x.
    \end{cases}
  \end{equation*}
  Let $\widecheck{\mathcal C}$ be a chain of subsets of $X$ such that
  for $Y\in \widecheck{\mathcal C}$, $k|_Y$ is uniformly regular.  If
  it were the case that $\bigcup_{Y\in\widecheck{\mathcal C}} Y = X$,
  then for some $Y$, $X\backslash Y$ would be countable, and hence $Y$
  would be dense in $[0,1]$.  But then we could find $x,y\in Y$ such
  that $1/|y-x|$ is arbitrarily large, contradicting the uniform
  regularity of $Y$.
\end{example}

\end{document}